\documentclass[12pt]{amsart}
\usepackage{amssymb,amsmath,amsfonts,amsthm,nomencl,mathrsfs} 
\usepackage[arrow, matrix, curve]{xy}
\usepackage[latin1]{inputenc}
\usepackage{a4wide}

\newcommand{\IR}{\mathbb{R}}

\newcommand{\question}[1]{\leavevmode{\marginpar{\tiny%
$\hbox to 0mm{\hspace*{-0.5mm}$\leftarrow$\hss}%
\vcenter{\vrule depth 0.1mm height 0.1mm width \the\marginparwidth}%
\hbox to 0mm{\hss$\rightarrow$\hspace*{-0.5mm}}$\\\relax\raggedright #1}}}

\newcommand{\ICC}{C^{\infty}}

\newcommand{\loc}{\mathrm{loc}}

\renewcommand{\c}{ c }

\newcommand{\IL}{L}

\newcommand{\Id}{d}

\newcommand{\f}{\frac}
\newcommand{\nn}{\nonumber}

\theoremstyle{plain}            
\newtheorem{theorem}{theorem}[section]

\newtheorem{Theorem}[theorem]{Theorem}
\newtheorem{Proposition}[theorem]{Proposition}

\newtheorem{Conjecture}[theorem]{Conjecture}

\theoremstyle{definition}       
\newtheorem{Definition}[theorem]{Definition}

\newtheorem{Remark}[theorem]{Remark}

\begin{document}

\setcounter{section}{1}

\begin{titlepage}
\title[]{The BMS conjecture }

\author[B. G\"uneysu]{Batu G\"uneysu}

\end{titlepage}

\maketitle 

\begin{abstract}
I explain an open conjecture by Braverman/Milatovic/Shubin (BMS) on the positivity of square integrable solutions $f$ of $(-\Delta+1)f\geq 0$ on a geodescially complete Riemannian manifold, and its connection to essential self-adjointness problems of covariant Schrödinger operators. The latter conjecture has remained open for more than 14 years now.
\end{abstract}

Let $M$ be a smooth connected Riemannian manifold, equipped with its usual Riemannian volume measure $\Id\mu$. We denote the scalar Laplace-Beltrami operator with\footnote{A \lq\lq{}$\dagger$\rq\rq{} always stands for the formal adjoint of a differential operator acting between sections of metric vector bundles over $M$; it depends on the fixed Riemannian metric on $M$ and the underlying metrics on the bundles (which are trivial in the scalar case).} $-\Delta=\Id^{\dagger}\Id$, and its Friedrichs realization in $\IL^2(M)$ with $H\geq 0$, where we understand all our spaces of functions to be complex-valued, unless otherwise stated. Let $E\to M$ be a smooth complex metric vector bundle with a smooth metric covariant derivative $\nabla$ thereon. Spaces of sections having a certain global or regularty $*$ will be denoted with $\Gamma_*$. For any Borel section $f$ of $E\to M$, the section
$\mathrm{sign}(f)\in\Gamma_{L^{\infty}}(M,E)$ is defined by 
$$
\mathrm{sign}(f)(x):=\begin{cases}&\f{f(x)}{|f(x)|}, \>\text{ if } \>f(x)\ne 0\\
&0, \>\text{ else.}\end{cases}
$$
The central result behind anything that follows is the following geometric variant of a classical distributional inequality by Kato \cite{kato2} (\lq\lq{}covariant Kato inequality\rq\rq{}): \emph{For all $f\in\Gamma_{\IL^1_{\loc}}(M,E)$ with $\nabla^{\dagger}\nabla f\in\Gamma_{\IL^1_{\loc}}(M,E)$ weakly, one has the weak inequality }
\begin{align}\label{fkqq}
-\Delta |f|\leq \Re\left( \nabla^{\dagger}\nabla f,\mathrm{sign}(f)\right).
\end{align}
A proof of the latter inequality can be found in \cite{braver}. It is in fact a local result which therefore holds without any further assumptions on $M$. Let us now pick a potential 
\begin{align}\label{poteea}
0\leq V\in\Gamma_{\IL^2_{\loc}}(M,\mathrm{End}(E))
\end{align}
 and assume we want to prove that the symmetric nonnegative operator $(\nabla^{\dagger}\nabla+V)|_{\Gamma_{\ICC_{\c}}(M,E)}$ in $\Gamma_{\IL^2}(M,E)$ is essentially self-adjoint. By an abstract functional analytic fact and some simple distribution theory, the latter essential self-adjoint is equivalent to the following implication:
\begin{align}\label{vor}
&f\in\Gamma_{\IL^2}(M,E),\>(\nabla^{\dagger}\nabla+V+1)f=0\>\text{ weakly }\\\nn
&\Rightarrow \> f=0.
\end{align}
So let $f$ be given with (\ref{vor}). In order to prove $f=0$, following Kato's original approach for $M=\IR^m$, it is tempting to use the covariant Kato inequality, which in combination with $V\geq 0$ immediately implies
$$
(-\Delta +1)(-|f|)\geq 0 \quad\text{weakly}.
$$

This motivates the $\mathscr{C}=L^2_{\IR}(M)$ case of following definition, which is taken from \cite{g54}:

\begin{Definition}\label{kl} Let\footnote{$L^1_{\mathrm{loc},\IR}$ stands for the space of \emph{real-valued} locally integrable functions, and likewise for $L^2_{\IR}(M)$.} $\mathscr{C}\subset L^1_{\mathrm{loc},\IR}(M)$ be an arbitrary subset. Then the Riemannian manifold $M$ is called $\mathscr{C}$-\emph{positivity preserving} (PP), if the following implication of weak inequalities holds true for every $\phi\in\mathscr{C}$,
\begin{align}\label{ppp}
 (-\Delta+1)\phi \geq 0 \Rightarrow  \phi\geq 0.
\end{align}
\end{Definition} 

Assume now $M$ is $L^2_{\IR}(M)$-positivity preserving. Then in the above situation we can conclude $-|f|\geq 0$, thus $f=0$, and we have shown:

\begin{Proposition} If $M$ is $L^2_{\IR}(M)$-positivity preserving, then for every potential $V$ with (\ref{poteea}), the operator $(\nabla^{\dagger}\nabla+V)|_{\Gamma_{\ICC_{\c}}(M,E)}$ in $\Gamma_{\IL^2}(M,E)$ is essentially self-adjoint.
\end{Proposition}

On the other hand, either using refined integration by parts techniques \cite{braver} or using wave equation techniques \cite{gerw}, one can prove:

\begin{Theorem} If $M$ is geodesically complete, then for every potential every potential $V$ with (\ref{poteea}), the operator $(\nabla^{\dagger}\nabla+V)|_{\Gamma_{\ICC_{\c}}(M,E)}$ in $\Gamma_{\IL^2}(M,E)$ is essentially self-adjoint.
\end{Theorem} 

This lead M. Braverman, O. Milatovic and M. Shubin to the following conjecture from 2002, which I formulate for convenience in the language of Definition \ref{kl}\footnote{Note that the BMS-conjecture is much older than Definition \ref{kl}, which was in fact modelled on the conjecture.}:

\begin{Conjecture}[BMS-conjecture] If $M$ is geodesically complete, then $M$ is $L^2_{\IR}(M)$-PP.
\end{Conjecture}

I invite the interested reader to attack this problem, which is still open in this generality!\vspace{2mm}

It is instructive in this context to explain Kato's simple and elegant proof of the fact that the Euclidean $M=\IR^m$ is $L^2_{\IR}(\IR^m)$-PP: In this case, $\Delta+1$ induces an isomorphism (of topological linear spaces)
$$
\Delta+1: \mathscr{S}(\IR^m)\rq{}\xrightarrow{\sim} \mathscr{S}(\IR^m)\rq{}
$$
on the space of Schwartz distributions, whose inverse is positivity preserving. Thus, if a real-valued $\phi\in\IL^2(\IR^m)\subset \mathscr{S}(\IR^m)\rq{}$ satisfies (\ref{ppp}), then we can immeadiately conclude $\phi\geq 0$.\vspace{2mm}

On a general Riemannian manifold there seems to be no appropriate substitute for the space of Schwartz distributions, and so one needs a new idea. The best result known so far on general Riemannian manifolds on the full $\IL^q$-scale is the following one from \cite{guenb} (which slightly generalizes \cite{g54}), that requires an additional lower bound on the Ricci curvature: 

\begin{Theorem}\label{cuttta} If $M$ is geodesically complete with a Ricci curvature bounded from below by a constant, then $M$ is $L^q_{\IR}(M)$-PP for all $q\in [1,\infty]$.
\end{Theorem}

The reader may find the following final remarks helpful:

\begin{Remark}1. The proof of Theorem \ref{cuttta} is based on the construction of a sequence of \emph{Laplacian cut-off functions} (cf. \cite{g54} for a precise definition), which leads to the assumption on the Ricci curvature. Once one has such a sequence, at least the $q=2$ case follows easily using that $(H+1)^{-1}$ is positivity preserving on $\IL^2(M)$ in combination with simple integration by parts arguments. The $q\neq 2$ case requires an additional argument to prove the boundedness of $\Id  (H+1)^{-1}$ from $\IL^q(M)$ to $\Omega^1_{\IL^q}(M)$, which again leads to the curvature assumption. If $M$ is the Euclidean $\IR^m$, such a sequence of Laplacian cut-off functions is readily obtained using the distance function and scaling.\\ 
2. The $q=2$ case from Theorem \ref{cuttta} can be generalized to allow a Ricci curvature having an appropriate \emph{variable} lower bound, as then one can still prove the existence of a sequence of Laplacian cut-off functions (cf. \cite{alberto}).\\
3. It really makes sense to consider the positivity preservation property on a full $\IL^q$-scale: For example, it is easy to check \cite{g54} that every $(\ICC\cap\IL^{\infty}_{\IR})$-PP Riemannian manifold is stochastically complete, meaning that
$$
\int_M \mathrm{e}^{-t H}(x,y) \Id\mu(y)=1\quad\text{ for all $t>0$, $x\in M$,}
$$
or in other words, that Brownian motions on $M$ cannot explode in a finite time. So for example, Theorem \ref{cuttta} provides an independent proof of S.T. Yau\rq{}s classical result which states that geodesically complete Riemannian manifolds with a Ricci curvature bounded from below by a constant are stochastically complete. This was my original motivation for the general form of Definition \ref{kl}, that is, the definition should be flexible enough to deal with problems such as stochastic completeness and essential self-adjointness simultaniously. 
\end{Remark}

\vspace{2mm}

{\bf Acknowledgements:} I would like to thank W. Arendt for his interest in my point of view on the BMS conjecture, and Alberto Setti for a very helpful correspondence on cut-off functions. The author research been financially supported by the DFG 647: Raum-Zeit-Materie.

\end{document}